\documentclass[12pt]{amsart}

\usepackage{amsgen,amsthm,amstext,amsbsy,amsopn,amsfonts,amssymb, enumerate}

\usepackage{amsmath}

\newcommand\la{\langle}
\newcommand\ra{\rangle}

\newcommand\nn{{\mathfrak n}}

\newcommand\vv{{\mathfrak v}}

\newcommand\zz{{\mathfrak z}}

\newcommand\NN{\mathbb N}

\newcommand\RR{\mathbb R}
\newcommand\ZZ{\mathbb Z}

\newcommand\ad{\operatorname{ad}}

\newcommand\grad{\operatorname{grad}}

\theoremstyle{plain}

\newtheorem{thm}{Theorem}[section]
\newtheorem{lem}[thm]{Lemma}
\newtheorem{prop}[thm]{Proposition}
\newtheorem{cor}[thm]{Corollary}

\theoremstyle{definition}
\newtheorem{defn}[thm]{Definition}
\newtheorem{rem}[thm]{Remark}
\newtheorem{example}[thm]{Example}

\begin{document}

\title[The  geodesic flow on nilmanifolds associated to graphs]
{The  geodesic flow on nilmanifolds associated to graphs}

\author{Gabriela P. Ovando}

\thanks{{\it (2000) Mathematics Subject Classification}: 70G65, 70H05, 70H06, 22E70, 22E25 }


\thanks{Work partially supported by ANPCyT, CONICET  and SCyT (UNR)}

\address{CONICET- UNR, Departamento de Matem\'atica, ECEN - FCEIA, Pellegrini 250, 2000 Rosario, Santa Fe, Argentina.}

\

\email{gabriela@fceia.unr.edu.ar}

\begin{abstract} In this work we study the geodesic flow on nilmanifolds associated to graphs.  We are interested in the construction of  first integrals to show  complete integrability on some compact quotients. Also  examples of integrable geodesic flows and of non-integrable ones are shown.    
\end{abstract}

\maketitle

 \noindent\section{Introduction}

In this work we deal with the integrability of the geodesic flow on nilmanifolds.  The condition of integrability   imposes  obstructions to the topology of the supporting manifold \cite{Ta1,Ta2}. In known examples constructed by starting with a Lie group $N$,  one is able to construct a set of linear independent first integrals in involution on the tangent space $TN$, and doing some extra work,  one can induce the first integrals to $T(\Gamma \backslash N)$, for a cocompact lattice $\Gamma < N$. These first integrals cannot be analytic in most of these examples. And in any case, the work to find enough functions in involution is doing by hand case by case, that shows that  the   integrability question is not well understood. Even under restricted conditions, for instance on locally homogeneous manifolds, there is no general theory to decide the complete integrability of the geodesic flow. In \cite{KOR} the authors studied some relationships between the algebra of first integrals and the algebraic geometric structure of the underlying Lie group, more specifically  the corresponding isometry group. 

Lie groups and their compact quotients were already used to answer nice geometrical questions \cite{BJ, BJ2, BT1, BT2,GMS, Pa,Sc}. 
For 2-step nilpotent  Lie groups Butler introduced the notion of non-integrable Lie algebras, and proved that the geodesic flow  on compact quotients manifolds arising from Lie groups with non-integrable Lie algebras cannot be integrable. It is not hard to see that non-integrable Lie algebras are singular. In the same work Butler \cite{Bu1} proved  non-commutative integrability for manifolds associated to almost non-singular Lie algebras. Since 
 Bolsinov and Jovanovic  \cite{BJ}  proved that integrability in the non-commutative sense implies Liouville integrability, one gets  the integrability of the geodesic flow on manifolds $T(\Gamma\backslash N)$, where  $N$  is associated to an almost non-singular Lie algebra. Thus we have some questions between the type of the Lie algebra $\nn$, as non-singular, almost non-singular or singular, and the integrability of the geodesic flow on $M:=T(\Gamma\backslash N)$:
\begin{center}
$\nn$ almost non-singular \quad $\Longrightarrow$ \quad integrability on $M$;
\end{center}
\begin{center} 
$\nn$ non-integrable \quad $\Longrightarrow$ \quad the geodesic flow on $M$ non-integrable. 
\end{center}
\begin{center} 
$\nn$ non-integrable \quad $\Longrightarrow$ \quad $\nn$ singular. 
\end{center}
In this work we study the integrability question on a family of nilmanifolds arising from graphs. This family was introduced by Dani and Mainkar in \cite{DM}  to study  Anosov automorphisms on nilmanifolds. More recently, in \cite{DDM} the authors find Heisenberg like algebras in this family and they  give conditions  on the graph $G$ and on a lattice $\Gamma\subset N$ for which the compact quotient $\Gamma \backslash M$,  has a dense set of smoothly closed geodesics. 

We start with preliminaries to study the geodesic flow on nilmanifolds. The  structure of $N$ is used when indentifying $TN$ with $T^*N$ via the metric.  More details can be found in \cite{KOR}. Summarising  we get results in the following:
\begin{itemize}
{ 
\item there is a family of 
compact manifolds $M=\Gamma\backslash N$ with integrable geodesic flow such that the Lie algebra of $N$ is singular. They are associated to the star graphs on $k+1$-vertices, $S_k$, for $k\geq 3$. For $k=2$ the Lie algebra arising from the graph is the Heisenberg Lie algebra of dimension three.  Topologically the compact quotients are $S^1$-fiber bundles over $T^{2k}$. 

\item Nilmanifolds $\Gamma \backslash N$ for $N$ constructed from complete graphs on $s$-vertices, $K_s$, are
\begin{itemize}
\item non-integrable for $s=2n+1$.  This generalizes the example in  \cite{Bu1}.
\item almost non-singular for $s=2n$. 
\end{itemize}

\item Let $G$ denote a connected graph on $k$ vertices with $k\leq 4$. Then except for the complete graph $K_3$, any 2-step nilpotent Lie group $N_G$ so as the corresponding compact quotient admits  a completely integrable geodesic flow. 

\item For the graph $G=K_3$, we prove that the geodesic flow on the corresponding manifold $TN_G$ is completely integrable, although it cannot be completely integrable on $T(\Gamma \backslash N_G)$ for any discrete cocompact subgroup $\Gamma < N_G$ \cite{Bu1,Bu3}. We find five linear independent invariant functions. 

\item For the path graph in four vertices $P$, although  integrability of the geodesic flow can be derived from the almost singularity property, we explicitly show  a set of first integrals in involution on $TN_P$. 
}
\end{itemize}

Tools and techniques used here have no direct relationship with the previous theory developped in the 80's. We refer to Hamiltonian sytems constructed with an algebraic data as in the Adler-Kostant-Symes scheme \cite{Ko2,Sy}, or the work by  Thimm concerning the geodesic flow, see \cite{Th}.
In the present paper first integrals are constructed either related to Killing vectors or as invariant functions. The invariant notion is attached to the  natural action of the Lie group $N$ on $TN$ induced by translations on the left by elements of the group. Notice that invariant functions are trivially induced to any tangent space $T(\Gamma\backslash N)$. The geometry of 2-step nilpotent Lie groups equipped with a left-invariant metric as well as  the geometry around closed geodesics on compact quotients was extensively studied in \cite{Eb}. 

In most examples given here, we show explicit computations, in the framework of Lie groups and Lie algebras. All  these objects are assumed over the real numbers.

 \section{Preliminaries}
In this section we recall basic notions for the study of the geodesic flow and we introduce the tools to study the geometry 2-step nilpotent Lie groups when equipped with a left-invariant metric. In fact this is determined at the Lie algebra level. We are interested in Lie algebras which can be constructed from graphs.

\subsection{The geodesic flow on nilmanifolds}

Here we consider the geodesic flow on Lie groups equipped with a left-invariant metric. We recall the general setting. 

Let $TN$ denote the tangent bundle of the Lie group $N$. The geodesic field is related to the Hamiltonian vector field of the energy function $E: TN \to \mathbb R$ (see for instance \cite{Eb2}), 
$$E(p, Y) = \frac{1}{2}\langle Y, Y\rangle.$$

 We identify $TN$ with $N\times\nn$. In fact for $p\in N$ and  $v_p\in T_pN$ we associate the pair $(p, Y)\in N\times \nn$ where $Y$ is the left-invariant vector field on $N$ such that $Y(p)=v_p$. Therefore making use of these identifications  one also has 
$$
T_{(p,Y)}(TN) \simeq \nn \times \mathfrak n  = \{(U, V): U, V \in \mathfrak n\}.
$$
The manifold $TN$ is considered with the product metric
$$\la (U,V), (U',V')\ra_{(p,Y)}=\la U,U'\ra+\la V, V'\ra$$
and it  has a canonical symplectic structure which is induced by the canonical symplectic structure on $T^*N$ via the metric. This gives the definition of a Poisson bracket on $C^{\infty}(TN)$ denoted by $\{\,,\,\}$ as 
$$\{f, g\} = \Omega(X_f, X_g), \mbox{ for } f,g\in C^{\infty} TN,$$
where $\Omega$ is the  canonical symplectic form  on $TN$,
and $X_f$, $X_g$ are the {\it Hamiltonian} vector fields of $f$ and $g$ respectively.
Recall that for a smooth function $h: TN \to \mathbb R$, the Hamiltonian vector field of $h$, denoted by $X_h$ is  implicitly  given by 
$$dh_{(p, Y)}(U, V) = \Omega_{(p, Y)}(X_h(p, Y), (U, V)).$$
On the other hand, 
the gradient field of $h$, denoted by $\grad h$ is given by the formula
$$dh_{(p, Y)}(U, V) = \langle \grad_{(p, Y)}h, (U,V)\rangle.$$

The geodesic field on $TN$ is the vector field associated with the geodesic flow
$$\Phi_t(p, Y) = \gamma'(t)$$
where $\gamma(t)$ is the geodesic on $N$ with initial conditions $\gamma(0) = p$, $\gamma'(0) = (p, Y)$.

We say that
 a smooth function $f: TN \to \mathbb R$ is a \emph{first integral} of the geodesic flow if $f$  is constant along the integral curves of the geodesic field, $X_E(f)=0$,  equivalently if 
$$\{f, E\} = 0.$$ 

Notice that the gradient field of the energy function for a left-invariant metric is 
$\grad E(p,Y)=(0,Y)$. 
The proof of the following proposition arises from the definitions above. See \cite{KOR}. 

\begin{prop} Let $N$ denote a   Lie group equipped with a left-invariant metric. Let $f, g \in C^\infty(TN)$ be smooth functions with $\grad_{(p, Y)}f = (U, V)$ and $\grad_{(p, Y)} g = (U', V')$. Then
\begin{enumerate}[(i)]
\item  the Hamiltonian vector field for a smooth $f: TN \to \RR$  is 
\begin{equation}\label{hamiltonian2}
  X_f(p, Y) = (V, \ad^t(V )(Y) - U),
\end{equation}
where $\ad^t(V)$ denotes the transpose of $\ad(V)$ with respect to the metric on $\nn$. 
\item  the Poisson bracket follows
\begin{equation}\label{poisson2}
\{f, g\}(p, Y) = \langle U, V'\rangle - \langle V, U'\rangle + \langle Y, [V', V] \rangle.
\end{equation}
 
\end{enumerate}
\end{prop}

Let $f, g: TN \to \RR$ be smooth functions. We say that they are {\em in involution} whenever they Poisson commute: 
$$\{f,g\}=0.$$

The definition above says that  $\{f,g\}=0$ if and only if $df(X_g)=0$  if and only if $X_g(f)=0$ (and also $X_f(g)=0$). Therefore $f$ is constant along integral curves of $X_g$ (analogously for $X_f$).

The question of explicitly finding   functions  in involution is a research topic with several open questions. 
If $M$ is a Riemannian manifold and $X^*$ is a Killing vector field on $M$, then the function $f_{X^*}: TM \to \mathbb R$ defined as $f_{X^*}(v) = \langle X^*(\pi(v)), v\rangle$ is a first integral of the geodesic flow. But in general it is not clear if one can produce enough functions in involution for proving the complete integrability. 

We shall discuss later the question of {\em invariant functions}. Indeed a Lie group $N$ acts on its tangent bundle $q \cdot v_p=dL_{q}v_p$, that in terms of the identification above gives
$$ q \cdot (p, Y)=(qp, Y).$$
A function $f:TN \to \RR$ is called {\em invariant} if $f(qp, Y)=f(p,Y)$ for all $p,q\in N$, that is, $f$ is invariant under the action of $N$.
 
\begin{defn}  Let ($M,\la\,,\ra)$ be a Riemannian manifold. It   has {\em completely integrable} geodesic flow (in the sense of Liouville) if there exist $n$ first integrals of the geodesic flow, $f_i: TM \to \RR$, where $n = \dim M$, such that $\{f_i, f_j\} = 0$ for all $i, j$ and the gradients of $f_1, \ldots, f_n$ are linear independent on an open dense subset of $TM$.
\end{defn}

Integrability or non-integrability are difficult questions. Integrability imposes topological restrictions on compact manifolds.  

Let $N$ denote a Lie group equipped with a left invariant metric $g$. We say that $\Lambda \subset N$ is a {\em lattice} if $\Lambda$ is a discrete subgroup such that the quotient $\Lambda \backslash N$ is a compact space. Sometimes one also says that $\Lambda$ is a discrete cocompact subgroup of $N$. 

Mal'cev \cite{Mal} has shown that a simply connected nilpotent Lie group $N$ admits a lattice if and only if its Lie algebra admits a basis with rational structure constants. 
In this case the metric $g$ is induced to the quotient and it is also denoted by $g$. It satisfies
$$g(dp_m X,dp_m Y)= g(X,Y) \qquad \mbox{ for } X, Y \in \nn, m\in N,$$
where $p: N \to \Lambda \backslash N$ is the canonical projection.
One looks for first integrals not only on $N$ but on the compact quotients $\Lambda \backslash N$.

\begin{thm}  (\cite{Bu3}, Theorem 1.3) Let $\nn$ be a non-integrable 2-step nilpotent Lie
algebra with associated simply connected Lie group $N$. Assume that there exists a discrete, cocompact subgroup $\Lambda$ of $N$. Then for any such $\Lambda$ and any
left-invariant metric $g$ on $N$, the geodesic 
flow of $(\Lambda \backslash N, g)$ is not completely
integrable.
\end{thm}

The precise definition of non-integrable 2-step nilpotent Lie algebra will be given in the next subsection. 

\subsection{2-step nilpotent Lie groups and graphs}

Let $N$ denote a $2$-step nilpotent Lie group equipped with a left-invariant metric $\la\,,\,\ra$. Its Lie algebra $\nn$  decomposes as the orthogonal sum
$$\mathfrak n = \mathfrak v \oplus \mathfrak z, \qquad \quad \mbox{where} \quad \vv=\zz^{\perp}$$
 where $\mathfrak z$ is the center of $\mathfrak n$. In this situation each element $Z \in \mathfrak z$ induces a skew-symmetric linear map on $\vv$,  $j(Z): \mathfrak v \to \mathfrak v$, given by
\begin{equation}\label{br}
  \langle j(Z)U, V\rangle = \langle[U, V], Z\rangle
\end{equation}
for all $U, V \in \mathfrak v$. The geometry of $N$ is encoded in the maps $j(Z)$ \cite{Eb}.

Recall that whenever $N$ is simply connected, the exponential map, $\exp:\nn\to N$, is a diffeomorphism with inverse map $\log: N \to \nn$.  Moreover one has the formula:
$$\exp(X) \exp(Y)= \exp( X+Y + \frac12 [X,Y]) \qquad \mbox{ for all } X, Y\in \nn. $$
This formula enables the realization of the Lie group at the Lie algebra level. In fact, via the exponential map we define a product on the Lie algebra:   take the left-invariant vector fields $X, Y \in\nn$ and define a product $X \cdot Y$ by $X+Y + \frac12 [X,Y]$. In this way  the exponential map, as the only one-parameter group on $N$ with initial condition $X\in \nn$, becomes the map that, whenever $X=\sum_i x_i X_i$, for a basis of vector fields $\{X_i\}$,  sends $X$   to the vector $(x_1, x_2, \hdots , x_n)$ in $\RR^n$, in  usual coordinates.

In particular, for a 2-step nilpotent Lie group, after Equations (\ref{hamiltonian2}) and (\ref{poisson2}), one has:
\begin{itemize}
\item For the energy function $E:TN \to \RR$, its Hamiltonian vector field is
$$X_E(p, Y) = (Y, j(Y_\zz)Y_\vv).$$
\item A function   $f: TN \to \RR$ with gradient $\grad f(p,Y)=(U,V)$ is a first integral of the geodesic flow on $TN$ if and only if 
\begin{equation}\label{E-geod}\la Y, U\rangle = \langle j(Y_{\zz}) V_{\vv}, Y_{\vv}\ra.
\end{equation}
\end{itemize}

A distinguished family of  2-step nilpotent Lie algebras can be constructed starting with a graph $G$. 
Let $G$ be a  directed graph with at least one edge. Denote the vertices of $G$ by $S=\{X_1, \hdots, X_m\}$,  the edges of $G$ 
by $E=\{Z_1, \hdots ,Z_q\}$. 

 The Lie algebra
$\nn_G$ is the vector space direct sum $\nn_G=\vv \oplus \zz$ 
where we let $E$ be a basis over $\RR$ for $\zz$ and $S$ be a basis over $\RR$ for $\vv$. Define the bracket relations among elements of $S$ according to adjacency rules:
\begin{itemize}
\item  if $Z_k$
is a directed edge from vertex
$X_i$ to vertex $X_l$ then define
the skew-symmetric bracket 
$[X_i, X_l] = Z_k.$
\item If there is no edge between two vertices, then
define the bracket of those two elements in
$S$ to be zero. 
\end{itemize}
Extend the bracket relation to
all of $\vv$ by using bilinearity of the bracket.  

\begin{rem} In \cite{Ma} it was proved that the 2-step nilpotent Lie algebras associated with two directed graphs are Lie isomorphic if and only if the graphs from which they arise are isomorphic.
\end{rem} 

Choose the inner product on $\nn_G$ so that $S\cup E$ is an orthonormal basis for
$\nn_G$.
Observe that if $Z_k$ is a directed edge from $X_i$ to $X_l$ then the map $j(Z_k)$ defined in Equation \ref{br} satisfies $j(Z_k)X_i=X_l$ 
and 
$j(Z_k) X_p = 0$ for any other $X_p\in S$ where $p\neq i, l$.

\begin{example} \label{star} 
The {\em star graph} $S_k$ has $k+1$ vertices $V_0, V_1, \hdots V_k$ and edges $Z_i$. The vertices $V_0$ and $V_i$ are joined by the  edge $Z_i$, so that the Lie bracket gives $[V_0,V_i]= Z_i$. 
Setting $Z=a_1 Z_1 + \hdots + a_k Z_k$, the matrix presentation of $j(Z):\vv \to \vv$ in the basis $V_0, V_1, \hdots, V_k$, is given by
$$\left( \begin{matrix}
0 & -a_1 & -a_2 & \hdots & -a_k \\
a_1 & 0 & 0 \hdots & 0 \\
a_2 & 0 & 0 & \hdots & 0 \\
\vdots & \vdots & \vdots & \ddots & \vdots\\
a_k & 0 & 0 & \hdots & 0 
\end{matrix}
\right).
$$ 
Therefore $j(Z)$ is singular for $k+1> 2,$ that is for the star graph  with $k+1> 2$ vertices. 
\end{example}

\begin{example} \label{complete} The {\em complete graph} on $n$ vertices $K_n$ is the graph which has an edge between every pair of distinct vertices.
The graph $K_2$ corresponds to the known Heisenberg Lie algebra of dimension three, it has two vertices $V_1, V_2$ and one edge $Z$, 
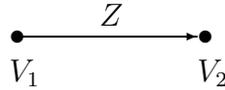
\begin{center}
\setlength{\unitlength}{5mm} 
\begin{picture}(17,3)
\put(6,1.5){$\bullet$}
\put(6,0.5){$V_1$}
\put(6.2,1.7){\vector(1,0){4.8}}
\put(11,1.5){$\bullet$}
\put(11,0.5){$V_2$}
\put(8.4,2){$Z$}
\end{picture} 
\end{center}
so that the corresponding map $j(Z)$ is  non-singular.
The complete graph $K_3$ has three vertices $V_1, V_2, V_3$ and three edges $Z_1, Z_2, Z_3$. 
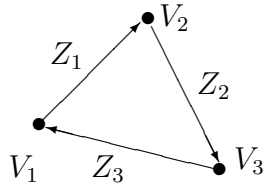
\begin{center}
\setlength{\unitlength}{5mm} 
\begin{picture}(17,5)
\put(6.1,1.7){$\bullet$}
\put(5.5,0.5){$V_1$}
\put(6.2,1.8){\vector(1,1){2.8}}
\put(6.6,3.5){$Z_1$}
\put(9,4.5){$\bullet$}
\put(9.5,4.5){$V_2$}
\put(10.5,2.7){$Z_2$}
\put(10.9,0.5){$\bullet$}
\put(11.5,0.7){$V_3$}
\put(7.7,0.5){$Z_3$}
\put(9.3,4.5){\vector(1,-2){1.8}}
\put(11,0.7){\vector(-4,1){4.5}}
\end{picture} 
\end{center}
Let $\nn_{K_3}$ denote the corresponding 2-step nilpotent Lie algebra where we have the Lie brackets
$$[V_1,V_2]=Z_1 \qquad [V_2,V_3]=Z_2 \qquad [V_3, V_1]=Z_3.$$
 Let $\la\,,\,\ra$ denote the metric on $\nn_{K_3}$ for which  this basis is orthonormal. The map $j(Z):\vv\to \vv$ has a matrix presentation in the basis $\{V_1,V_2,V_3\}$ of $\vv$ given by
$$j(aZ_1+bZ_2+cZ_3) = \left( \begin{matrix} 
0 & -a & -c \\
a & 0 & -b \\
c & b & 0 
\end{matrix}
\right).
$$ 
Notice that  the the dimension $ker j(Z)=1$,  for every $Z\in\zz-\{0\}$.
More generally consider the complete graph on $n$ vertices, $K_n$. The dimension of the center is $\dim \zz= \left( \begin{matrix} n \\
2 \end{matrix}\right)$ = $\dim \mathfrak{so}(n)$. Thus
\begin{itemize}
\item if $n$ is odd,  every $j(Z)$ is a singular map for every $Z\in \zz$; 
\item if $n$ is even, there exists $Z\in \zz$ such that $j(Z)$ is non-singular. In fact, assuming $V_1, V_2, \hdots V_{2s}$ are vertices and $Z_{ij}=[V_i, V_j]$, take the non-singular  map j(Z) with matrix
$$j(Z_{12}+Z_{34} + \hdots Z_{2s-1, 2s}) = \left(
\begin{matrix}
0 & -1 &  \\
1 & 0 & \\
 &  & 0 & -1 &\\
 &  & 1 & 0 & \\
 & & & & \ddots & \\
 & & & & &  &  0 & - 1\\
 & & & & &  & 1 & 0
\end{matrix}
\right).
$$
\end{itemize}

\end{example}

We say that a 2-step nilpotent Lie algebra $\nn$ is 

\begin{itemize}
\item  {\em non-singular} if $\ad(X) : \nn \to \zz$ is surjective for all $X \notin \zz$;
\item  {\em almost non-singular} if $j(Z)$  is non-singular for every $Z$ in an open dense subset of $\zz$.
\item {\em singular} if $j(Z)$ is singular for all $Z$ in $\zz$.
\end{itemize}

Every 2-step nilpotent Lie algebra belongs to one and only one of the types non-singular, almost non-singular
or singular \cite{GM}. 
 
 Notice that
\begin{itemize}
	\item the fact of being 
$\nn$  non-singular is equivalent to ask  $j(Z)$ to be non-singular for any $Z \in \zz -\{0\}$ for a (any) metric on $\nn$ (in fact this
does not depend on the choice of the left-invariant metric, see for instance \cite{Eb});
\item whenever a 2-step nilpotent Lie algebra $\nn$ is equipped with a metric, if there are
two nonzero elements $Z,Z'\in \zz$ such that $j(Z)$ is non-singular and $j(Z')$ is singular, then
the Lie algebra $\nn$ is almost non-singular.
\end{itemize}
\begin{rem}  Assume that $G$ is a graph with  at least one edge. If the graph $G$ is isomorphic to the complete graph $K_2$ then its Lie algebra $\nn_G$  is non-singular (this Lie algebra is also isomorphic to the Heisenberg Lie algebra of dimension three).

 Assume $G$ is not isomorphic to $K_2$. Thus $G$ contains an edge
$Z$ and a vertex $X$ for which the edge $Z$
is not incident to the vertex
$X$. Then $j(Z)X= 0$,  hence $\nn_G$ 
is either almost non-singular or singular. Lemma 3.3 in \cite{DDM} proved that the  Lie algebra $\nn_G$ is non-singular if
and only if $G = K_2$. 

Further if $G$ has more than one connected component, $\nn_G$ is almost non-singular if and only if each connected component is either non-singular or almost non-singular  \cite{DDM}.

\end{rem}

\begin{defn} \label{non-integrable} Let $\nn$ be a 2-step nilpotent Lie algebra, $\nn^*$ its dual space and for $\lambda\in \nn^*$ let $\ad^*(X)\lambda$ denote the element in $\nn^*$ given by $\ad^*(X)\lambda(Y)=-\lambda([X,Y])$, for all $X,Y\in \nn$.
\begin{enumerate}[(i)]
\item For $\lambda \in\nn^*$, let $\nn_{\lambda} := \{ X \in \nn\, :\,  \ad^*(X)\lambda = 0\}$.
\item  A $\lambda \in \nn^*$ 
is called {\em regular} if $\nn_{\lambda}$ has minimal dimension.
\item  A pair $\mu, \lambda\in\nn^*$ is called {\em generic} if $\dim [\nn_{\lambda}, \nn_{\mu}]$ is minimal. The Lie algebra $\nn$ is called {\em non-integrable} if for  a dense open subset of generic pairs $\lambda, \mu$, one has  $[\nn_{\lambda}, \nn_{\mu}]\neq 0.$ 
\end{enumerate}
\end{defn}

\begin{rem} Let $\lambda \in \nn^*$. Indeed   $\nn_{\lambda}$ is the isotropy algebra for the coadjoint representation: $p \cdot \lambda = -\lambda \circ Ad(g^{-1})$, for $p\in N$. Notice that $\lambda$ is regular whenever the dimension of the orbit of $\lambda$ under the coadjoint representation is maximal. It is clear that the center $\zz\subset \nn$ is contained in $\nn_{\lambda}$ for any $\lambda \in \nn^*$.
\end{rem}

We can read the non-integrability notion making use of tools at the Lie algebra level.  Let $\la\,,\,\ra$  denote an inner product on  the Lie algebra $\nn$, with respect to which there is the following  splitting as  direct sum of vector spaces $\nn=\vv \oplus \zz$, with $\vv=\zz^{\perp}$. 

For any $\lambda\in\nn^*$ there exists  unique $Z\in \zz,$ and $V\in \vv$ such that  $\lambda=\ell_{V+Z}$ where $\ell_{V+Z}(X)=\la V+Z, X\ra$. Thus we  denote $\nn_{\ell_{V+Z
}}$  directly by $\nn_{V+Z}$. So
$$\begin{array}{rcl}
\nn_{V+Z} & = & \{ X\in\nn\,:\, \la V+Z, \ad(X) U\ra=0 \mbox{ for all } U\in \nn\}\\
 & = & \left\{
\begin{array}{ll}
\nn & \mbox{ if } Z = 0\\
\zz \oplus ker\, j(Z) & \mbox{ if } Z\neq 0
\end{array} 
\right.
 \end{array}
$$

Assume $\ell_{V+Z}$ is regular,  thus 
for $\nn$ non-singular or almost non-singular, it is clear that $\nn_{V+Z}=\zz$. In both cases  for both $\lambda, \mu$ regular, one has $[\nn_{\lambda}, \nn_{\mu}]=0$. 

\begin{cor} \label{coro1} Let $\nn$ denote a non-integrable Lie algebra $\nn$ then $\nn$ must be singular. 
\end{cor}

\begin{example} In \cite{Bu1} Butler showed an example of a non-integrable Lie algebra. Its Lie group has Lie algebra isomorphic to $K_3$.
 
More generally assume that $n$ is odd,  $n=2k+1>2$, and $G=K_{n}$ is the complete graph. As mentioned above the Lie algebra $\nn_G$ is singular. Roughly speaking the graph $K_{2k+1}$ can be constructed from the complete graph $K_{2k}$ by adding one vertex and all the edges joining the added vertex with the previous ones.  

Let $S$ denote the set of vertices, assume $|S|=n$  with $n=2k+1$ and choose the vertex $V_{2k+1}\in S$. Denote by $Z_{ij}$(or $Z_{i,j}$) the basis element in $\zz$ such that $Z_{i,j}:=[V_i, V_j]$. 

Take $Z\in \zz$ defined as $Z:=\sum_{i=1}^k  Z_{2i-1, 2i}\in \zz$ Then the restriction of $j(Z)$  to the vector subspace $\mathfrak w_1$ spanned by $V_1, V_2, \hdots, V_{2k}$ - see matrix in Example \ref{complete} -, is a non-singular linear map. And  $V_{2k+1} \in ker j(Z)$ since  $j(Z_{2i-1, 2i})V_{2k+1}=0$ for all $i$. 

Analogously take $V_1$ in the kernel of $j(\tilde{Z})$ for $\tilde{Z}=\sum_{i=1}^k  Z_{2i, 2i+1}$ so that the restriction of $j(\tilde{Z}$) is non-singular on the vector subspace $\mathfrak w_2$ of dimension $2k$ spanned by $V_2, V_3, \hdots V_{2k+1}$. 

It is clear that $\ell_{Z}$ and $\ell_{\tilde{Z}}$ are regular in $\nn_G^*$. Moreover $V_1\in \nn_{\tilde{Z}}$, $V_{2k+1}\in \nn_Z$ and  $[V_1, V_{2k+1}]=Z_{1,2k+1}$. This implies that $K_{2k+1}$ is non-integrable. 
\end{example}

\begin{prop} \label{NonintegrableKn} Let $n=2k+1 \in \NN$ with $k\geq 1$ and let $G=K_n$ denote the complete graph on $n$ vertices. The corresponding Lie algebra $\nn_G$ is non-integrable.
\end{prop}

 Proposition above and Theorem 1.3 in \cite{Bu3} imply that the 2-step nilpotent Lie group $N_G$, constructed by starting with the complete graph $G=K_n$, with $n$ odd, has no completely integrable geodesic flow on $\Lambda \backslash N_G$, for any left-invariant metric and any lattice $\Lambda \subset N_G$.

\begin{example} \label{graphsandvertex}  Connected graphs with $n$  vertices and the type of their Lie algebras: 
\begin{enumerate}[(i)]
\item $n=2$: the complete $K_2$, $\nn_{K_2}$ non-singular,
\item $n=3$: the complete $K_3$, with $\nn_{K_3}$ singular, the star graph $S_3$ with $\nn_{S_3}$ singular;
\item $n=4$: The star graph $S_4$ with $\nn_{S_4}$ singular and the (corresponding to) almost-non-singular Lie algebras: the complete graph $K_4$, the cycle $C_4$, the path in four vertices $P_4$, and the graphs  $G_1$ and $G_2$. See \cite{DDM}. The graphs $G_1$ and $G_2$ are the two non-isomorphic graphs that one obtains as subgraphs of  $K_4$, by erasing  one or two edges.  
\end{enumerate}
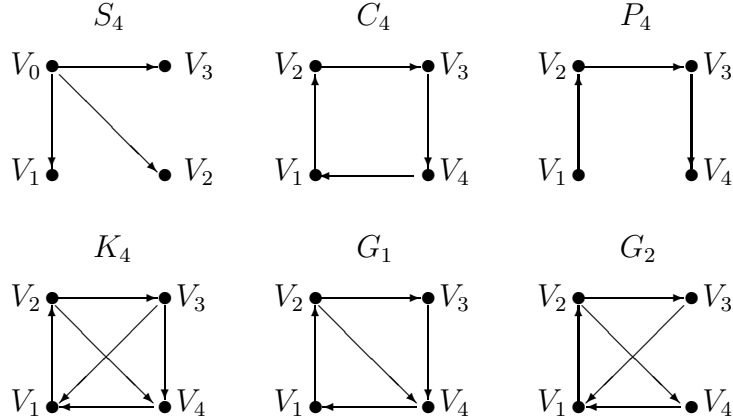
\begin{center}
\setlength{\unitlength}{5mm} 
\begin{picture}(17,5)
\put(2.3,4.3){$S_4$}
\put(1,0.1){$\bullet$}
\put(0.1,0.2){$V_1$}
\put(1.2,3){\vector(0,-1){2.5}}
\put(1,3){$\bullet$}
\put(0.1,3){$V_0$}
\put(4,3){$\bullet$}
\put(4.7,3){$V_3$}
\put(4,0.1){$\bullet$}
\put(4.7,0.2){$V_2$}
\put(1.2,3.2){\vector(1,0){2.8}}
\put(1.4,3){\vector(1,-1){2.6}}
\put(9.3,4.3){$C_4$}
\put(8,0.1){$\bullet$}
\put(7.1,0.2){$V_1$}
\put(8.2,0.3){\vector(0,1){2.7}}
\put(8,3){$\bullet$}
\put(7.1,3){$V_2$}
\put(11,3){$\bullet$}
\put(11.5,3){$V_3$}
\put(11,0.1){$\bullet$}
\put(11.5,0.2){$V_4$}
\put(8.2,3.2){\vector(1,0){2.8}}
\put(11.2,3){\vector(0,-1){2.5}}
\put(10.8,0.3){\vector(-1,0){2.5}}
\put(16.3,4.3){$P_4$}
\put(15,0.1){$\bullet$}
\put(14.1,0.2){$V_1$}
\put(15.2,0.3){\vector(0,1){2.7}}
\put(15,3){$\bullet$}
\put(14.1,3){$V_2$}
\put(18,3){$\bullet$}
\put(18.5,3){$V_3$}
\put(18,0.1){$\bullet$}
\put(18.5,0.2){$V_4$}
\put(15.2,3.2){\vector(1,0){2.8}}
\put(18.2,3){\vector(0,-1){2.5}}
\end{picture} 
\end{center}

\medskip

\begin{center}
\setlength{\unitlength}{5mm} 
\begin{picture}(17,5)
\put(2.3,4.3){$K_4$}
\put(1,0.1){$\bullet$}
\put(0.1,0.2){$V_1$}
\put(1.2,0.3){\vector(0,1){2.7}}
\put(1,3){$\bullet$}
\put(0.1,3){$V_2$}
\put(4,3){$\bullet$}
\put(4.5,3){$V_3$}
\put(4,0.1){$\bullet$}
\put(4.5,0.2){$V_4$}
\put(1.2,3.2){\vector(1,0){2.8}}
\put(4.2,3){\vector(0,-1){2.5}}
\put(3.9,0.3){\vector(-1,0){2.5}}
\put(1.3,3){\vector(1,-1){2.6}}
\put(4,3){\vector(-1,-1){2.6}}
\put(9.3,4.3){$G_1$}
\put(8,0.1){$\bullet$}
\put(7.1,0.2){$V_1$}
\put(8.2,0.3){\vector(0,1){2.7}}
\put(8,3){$\bullet$}
\put(7.1,3){$V_2$}
\put(11,3){$\bullet$}
\put(11.5,3){$V_3$}
\put(11,0.1){$\bullet$}
\put(11.5,0.2){$V_4$}
\put(8.2,3.2){\vector(1,0){2.8}}
\put(11.2,3){\vector(0,-1){2.5}}
\put(10.9,0.3){\vector(-1,0){2.5}}
\put(8.3,3){\vector(1,-1){2.6}}
\put(16.3,4.3){$G_2$}
\put(15,0.1){$\bullet$}
\put(14.1,0.2){$V_1$}
\put(15.2,0.3){\vector(0,1){2.7}}
\put(15,3){$\bullet$}
\put(14.1,3){$V_2$}
\put(18,3){$\bullet$}
\put(18.5,3){$V_3$}
\put(18,0.1){$\bullet$}
\put(18.5,0.2){$V_4$}
\put(15.2,3.2){\vector(1,0){2.8}}
\put(17.9,0.3){\vector(-1,0){2.5}}
\put(15.3,3){\vector(1,-1){2.6}}
\put(18,3){\vector(-1,-1){2.6}}
\end{picture}
\end{center}
\end{example}

\section{Involution of invariant functions}

The goal now is  the study of invariant functions.
Invariant functions   descend to any compact quotient $\Lambda\backslash N$ for any lattice $\Lambda\subset N$. Therefore it is desirable to have a good number of independent invariant functions. 

 Under the  natural action of $N$ on $TN\simeq N \times \nn$  given by $n \cdot (p,Y)= (np,Y)$, a function     $f:TN \to \RR$ is { invariant} if $f(p,Y)=f(e,Y)$ for all $p\in N, Y\in \nn$.
For instance if the metric on the Lie group $N$ is left-invariant, the corresponding energy function is invariant.

\label{finvariants} Concerning invariant functions in $C^{\infty}(TN)$ for a Lie group $N$,  from notions and properties above one proves   the next statements.
\begin{enumerate}[(i)]
\item The gradient field for an invariant function $f:TN \to  \RR$ has the form 
$$\grad_{(p,Y)}(f)=(0,V)$$ for some $V\in \nn$. In fact,  denote by $U,V$ the components of the gradient: $\grad f(p,Y)=(U,V)$. Since $f(e,Y)=f(p,Y)$ one has

$df_{(p,Y)}(U',V')=\frac{d}{ds}|_{s=0}f(e, Y+sV)=\la U,U'\ra + \la V,V'\ra$, so that $U=0$. 
The corresponding Hamiltonian vector field is given by 
$$ X_f(p,Y)=(V, \ad^t(V)Y)$$
where  $\ad^t(V)$ denotes the transpose of $\ad(V)$ relative to the metric on $\nn$.
 In particular an invariant function $f:TN \to \RR$ is a first integral of the geodesic flow if and only if 
$$0=\la Y, [V, Y]\ra \quad \mbox{ for } \quad (0,V)=\grad f(p,Y).$$ 

\item The set of invariant functions $\{f:TN \to \RR \, : \, f \mbox{ is invariant }\}$ is in correspondence with the set of functions on $\nn$: $\{ F: \nn \to \RR\}$.

Given an invariant function $f:TN \to \RR$ define $F: \nn \to \RR$ as
$F(Y)=f(e,Y)$
and conversely given $F:\nn\to \RR$ define an invariant function $f:TN  \to \RR$ by
$$f(p,Y)=F(Y)\qquad \quad \mbox{ for all } p\in N, Y\in  \nn.$$
\item  Let $f_1, f_2: TN \to \RR$ be invariant functions. Then their corresponding gradients follow $\grad(f_i)(p,Y)=(0,V_{F_i})$ for $i=1,2$ and the Poisson bracket is
\begin{equation}\label{bracket-invariants}\{f_1,f_2\}(p,Y)= - \la Y, [V_{F_1}, V_{F_2}]\ra, 
\end{equation}
where $V_{F_i}=\grad_{\nn}F_i$. 
 In fact on the one hand $\grad f_i(p,Y)=(0,V_i)$, so that 
$$\{f_1,f_2\}(p,Y)  =   - \la Y, [V_1, V_2]\ra, $$ 
but one also knows that 
$\frac{d}{ds}|_{s=0}f_i(e, Y+sV')= \frac{d}{ds}|_{s=0}F_i(Y+sV')= {dF_i}_{_Y}(V')=\la \grad_{\nn} F_i(Y), V'\ra$. These equations prove the Equality \ref{bracket-invariants}.

\end{enumerate}

The next proposition  specifies some invariant functions which are first integrals of the geodesic flow. Notice that for every $(p,Y)\in TN$, it holds $\grad f_{Z_0}(p,Y)=(0,Z_0)$ and $\grad g_A(p,Y)=(0,AY)$. See the proof  in \cite{KOR}. 

\begin{prop} \label{propcuad}
  Let $(N, \la\cdot,\cdot\ra)$ be a Lie group with a left-invariant metric.
  \begin{enumerate}[(i)]
  \item The function $f_{Z_0}: TN \to \mathbb R$, defined by $f_{Z_0}(p, Y) = \langle Y, Z_0\rangle$, is a first integral of the geodesic flow for all $Z_0 \in \mathfrak z$. Moreover, the family $\{f_{Z_0}\}_{Z_0 \in \zz}$ is a Poisson-commutative family of first integrals.
  \item Let $A: \mathfrak n \to \mathfrak n$ be a symmetric endomorphism of $\mathfrak n$ and let $g_A:TN \to \mathbb R$ denote the quadratic polynomial given by $g_A(p, Y) = \frac{1}{2}\langle Y, AY\rangle$. Then $g_A$ is a  first integral of the geodesic flow if and only if $\langle Y, [AY, Y]\rangle = 0$ for all $Y \in \nn$.
 \end{enumerate}
\end{prop}

Let $N$ denote a Lie group with Lie algebra $\nn$. The {\em Gauss map} $G:TN \to \nn$ is given by $G(p,Y)=Y$. Thus its  differential $dG_{(p,Y)}(U,V)=V$. Geometrically the Gauss map sends an element $Y_p\in TN$ to $dL_{p^{-1}}Y_p$ and takes this as the initial value corresponding to the left-invariant vector field $Y$. 
A smooth function $F:\nn\to \RR$, corresponds to a  smooth function $f:TN \to \RR$ defined as $f=F\circ G$.  Thus for $F_1, F_2\in C^{\infty}(\nn)$
$$\begin{array}{rcl}
\{F_1 \circ G, F_2\circ G\}(p,Y) & = & \{F_1, F_2\} \circ G (p,Y)\\
& = & -\la Y, [V_{F_1}, V_{F_2}]\ra.
\end{array}
$$

 The Gauss map in this context was previously studied en \cite{Eb2}. We notice that the Poisson bracket on the Lie algebra $\nn$ has symplectic leaves given by the ``coadjoint'' orbits, which are induced from $\nn^*$ to $\nn$ via the corresponding metric. However the Gauss map $G$ does not send $(p,Y)\in TN$ to  a vector in $\nn$ which is tangent to a coadjoint orbit, that is, one should project it if necessary.

For invariant functions, the Poisson bracket on $\nn$ (see \cite{Eb}) is in correspondence with the Poisson bracket on $TN$, $\{f_1,f_2\}(p,Y)= \{F_1,F_2\}(e,Y)$, 
and via the relation $f_i(p,Y)=F_i(Y)=F_i \circ G(p,Y)$, one has $\{f_1,f_2\}(p,Y)=\{F_1, F_2\}\circ G(p,Y).$ We summarize below these results. 

\begin{prop} Let $N$ denote a Lie group with Lie algebra $\nn$. Let $G:TN \to \nn$ denote the Gauss map. 
\begin{enumerate}[(i)]
\item For functions $F_1,F_2:\nn \to \RR$ one has $F_i\circ G :TN \to \RR$ so that for every $(p,Y)\in TN$: 
\begin{equation}\label{Poisson-invariant}
\{F_1 \circ G, F_2\circ G\}(p,Y)  =  \{F_1, F_2\} \circ G (p,Y) = -\la Y, [V_{F_1}, V_{F_2}]\ra
\end{equation}
where $V_{F_i}=\grad_{\nn} F_i$ so that $\grad(F_i\circ G)=(0, V_{F_i})$. Clearly $F_i\circ G: TN \to \RR$ is invariant by construction. 
\item Given a smooth function $f:TN \to \RR$ there exists $F:\nn \to \RR$ so that $F\circ G=f$ if and only if $f$ is invariant. In this case we have the commuting diagram
\begin{center}
\setlength{\unitlength}{5mm} 
\begin{picture}(17,7)
\put(6,5.5){$TN$}
\put(7.6,5.7){\vector(1,0){2.5}}
\put(8.5,6){$G$}
\put(11,5.5){$\nn$}
\put(5.2,3.6){$f$}
\put(6.5,5.2){\vector(0,-1){2.3}}
\put(10.2,5){\vector(-1,-1){2.7}}
\put(6.3,1.6){$\mathbb R$}
\put(8.3,3.7){F}
\end{picture} 
\end{center}
 So for   smooth invariant functions  $f_i:TN \to \RR$, with associated functions $F_i:\nn \to \RR$: where $f_i=F_i\circ G$, for $i=1,2$, the Poisson bracket follows as in Equation (\ref{Poisson-invariant}): $\{f_1,f_2\}(p,Y)=\{F_1, F_2\}\circ G(p,Y).$
\end{enumerate}
\end{prop}

In the following paragraphs we shall show a construction of invariant functions for 2-step nilpotent Lie groups.  

Assume $N$ is a 2-step nilpotent Lie group with Lie algebra $\nn$. Let $\la\,,\,\ra$ denote a left-invariant metric on $N$ which gives the orthogonal decomposition 
$$\nn=\vv \oplus \zz \qquad \quad \mbox{ for } \vv=\zz^{\perp}.$$

Assume $A:\nn \to \nn$ is a symmetric linear map, preserving this decomposition. In particular $A\vv \subseteq  \vv$. Thus the map $g_A$ given by $g_A(p,Y)=\la AY,Y\ra$ is a first integral of the geodesic flow if and only if

\smallskip

$\la j(Y_{\zz}) A Y_{\vv}, Y_{\vv}\ra=0 \qquad \mbox{ for all } Y=Y_{\zz} + Y_{\vv}\in \nn$

\smallskip

which is equivalent to:

\smallskip

$j(Z) A = A j(Z)$ for all $Z\in \zz$. 

\smallskip

Moreover the maps $g_A, g_B$ are involution (for respective symmetric maps $A,B:\vv \to \vv$) if and only if

\smallskip

$j(Z) AB = j(Z)B A$ for all $Z \in \zz$. See \cite{KOR} for details. 

 \smallskip

L. Butler in \cite{Bu1}  exhibited the following polynomial first integrals for the geodesic flow. The statement given in \cite{Eb2} is more adequate for our work here. 

\begin{prop} \label{fin} Let $N$ denote an almost non-singular  2-step nilpotent Lie group equipped with a left-invariant metric $\la\,,\,\ra$. Let $\nn$ denote its Lie algebra with orthogonal splitting $\nn=\vv \oplus \zz$, where $\dim \vv= 2n$. For $i=1, \hdots, n$, the invariant functions $f_i:TN \to \RR$ given by
\begin{equation}
f_i(p,Y)=\la V, j(Z)^{2i} V\ra \qquad Y = V + Z\in \nn=\vv\oplus \zz;
\end{equation}
are first integrals of the geodesic flow. 
\end{prop}

\begin{example} \label{path} Let $P$ denote the graph which is the  path of length three on four distinct vertices. Assume the vertices are $V_1, V_2, V_3, V_4$ and let $\nn_{P}$ be the 2-step nilpotent Lie algebra associated with this graph with the Lie brackets
$$[V_1, V_2]=Z_1\qquad [V_2,V_3]= Z_2\qquad [V_3,V_4]=Z_3.$$
Assume  a metric that makes of the set $V_1,V_2,V_3,V_4,Z_1,Z_2,Z_3$  an orthonormal basis. 

Then $j(Z_i)$ is singular for every $i$ but $j(Z_1+Z_3)$ is non-singular. Thus $\nn_{P}$ is almost non-singular. Note that a general $j(aZ_1 + bZ_2+c Z_3)$ in the basis $V_1, V_2, V_3, V_4$
has a matrix of the form
$$j(Z)=\left(
\begin{matrix}
0 & -a & 0 & 0\\
a & 0 & -b & 0 \\
0 & b & 0 & -c\\
0 & 0 & c & 0 
\end{matrix}
\right)
$$
It is not hard to prove that a symmetric map on $\vv$ inducing an invariant function on $T N_P$ of the form $g_A(p,Y)=\frac12 \la AY, Y\ra$ gives a first integral of the geodesic flow  only   for $A=Id$. In fact one verifies that $Id$ is the only solution for symmetric maps $S:\vv \to \vv$ satisfying the condition $[S, J(Z)]=0$ for all $Z\in \zz$. Thus we consider 
 the function $\bar{h}:\nn_P \to \RR$ given by 
$\bar{h}(V+Z)= \la V, j(Z)^{2}V\ra$ 
where 
$$j(Z)^2 = \left(
\begin{matrix}
-a^2 & 0 & ab & 0\\
0 & -a^2-b^2 & 0 & bc \\
ab & 0 & -b^2-c^2 & 0\\
0 & bc & 0 & -c^2 
\end{matrix}
\right).
$$
This induces an invariant first integral on $T N_P$ given by $h(p,Y)=\bar{h}(Y)$ as in Proposition \ref{fin}. The gradient of $h$ is given by

$\grad h(p,Y)=(0, 2([j(Z) V,V] + j(Z)^2 V)$, where $Y=V+Z\in \nn_P$. Summarising we get  the following set of invariant functions in involution:
\begin{itemize}
\item $g(p,Y)=\frac12 \la Y, Y\ra$;
\item $h(p,Y)= \la Y, j(Z)^{2}Y_{\vv} \ra$, where $Y_{\vv}$ denotes the projection of $Y\in \nn_P$ onto $\vv$;
\item  the three functions for the center $f_i(p,Y)=\la Y, Z_i\ra$ for $i=1,2,3$.
\end{itemize}
 Invariant functions clearly descend to any space $T(\Lambda \backslash N_P)$ for any cocompact lattice $\Lambda$, since the Lie algebra has rational structure constants. We need two more functions to have the complete integrability (in the Liouville sense) in this case.  
\end{example}

\begin{rem} The coadjoint action of $N$ to $\nn^*$ is induced to $\nn$ via the metric obtaining:
$$p \cdot Y = Ad^t(p^{-1})(Y) \qquad \quad \mbox{ for all }g\in N, X\in \nn,$$
where $Ad^t(p)$ denotes the transpose of the Adjoint map relative to the metric, $\la Ad^t(p)(Y), Y\ra = \la Y, Ad(p) X\ra$ for all $p\in N, X, Y\in \nn$. 
Also tangent vectors to the orbit are induced by vector fields on $\nn$, that is $X\in \nn$ gives $\tilde{X}(Y)$:
$$\tilde{X} (Y)=\frac{d}{ds}_{|_{s=0}} \exp(sX) \cdot Y = -\ad^t(X) (Y).$$
Given a function $F:\nn \to \RR$ one can consider its restriction to a coadjoint orbit. But functions which are not trivial on $\nn$ can became constant on a coadjoint orbit. Take for instance $\nn$ a 2-step nilpotent Lie algebra equipped with a metric and consider the function $F:\nn \to \RR$ given by $F_0(Y)=\la Y,Z_0\ra$, where $Z_0 \in \zz$ is a fixed vector. It is not hard to see that $\grad_{\nn} F_0(Y)=Z_0$. But the  restriction of $F_0$ to the coadjoint orbit, gives
$${dF_0}_Y(\tilde{X}) = \frac{d}{ds}|_{s=0} \la Ad^t(\exp -sX)(Y), Z_0\ra = -\la Y, \ad(X)Z_0\ra =0$$
which proves the assertion. Notice that the example holds for any Lie algebra with non-trivial center. 
\end{rem}

\subsection{Once the situation of $K_3$} Here we come back to   the 2-step nilpotent Lie group and  quotients associated to the complete graph $K_3$. The corresponding  simply connected Lie group $N_{K_3}$, admits a lattice $\Gamma$ and as already said, no quotient $\Gamma \backslash N_{K_3}$ can be endowed with Riemannian metric - induced from a left invariant metric on $N_{K_3}$ - such that the corresponding geodesic flow is completely integrable. This question was discussed in \cite{Bu3}. However the geodesic flow on $N_{K_3}$ can be completely integrable. 

Take the notation of Example \ref{complete}. Choose the metric on the Lie algebra $\nn_{K_3}$ for which the set $V_1,V_2,V_3,Z_1,Z_2,Z_3$ is an orthonormal basis. By Proposition \ref{propcuad} one has  the  set of invariant functions on $TN_{K_3}$ given by

$\begin{array}{rcl}
2 E(p,Y) & = & \sum_{i=1}^3 \la Y, V_i\ra^2 + \sum_{i=1}^3 \la Y, Z_i\ra^2\\
f_{Z_j}(p,Y) & = & \la Z_j, Y\ra \quad \mbox{ for } j=1,2,3.
\end{array}
$

\smallskip

And the  following invariant function 
\begin{equation}\label{G}
G(p,Y)= \la Y, Z_1\ra \la Y, V_3\ra + \la Y, Z_2\ra \la Y, V_1\ra + \la Y, Z_3\ra \la Y, V_2\ra
\end{equation}
 is also  a first integral of the geodesic flow. The gradient of $G$ is given by 
$$\grad G(p,Y)= (0,\la Y,Z_2\ra V_1 + \la Y,Z_3\ra V_2 + \la Y,Z_1\ra V_3 + \la Y,V_3\ra Z_1 + \la Y,V_1\ra Z_2 + \la Y,V_2\ra Z_3),$$
so that it is not hard to prove that it satisfies Condition  \ref{E-geod}:
$$
\la Y, [Y, \la Y,Z_2\ra V_1 + \la Y,Z_3\ra V_2 + \la Y,Z_1\ra V_3]\ra=0.
$$
On the other hand it is clear that $G$ is in involution with $f_{Z_1}, f_{Z_2}, f_{Z_3}$. 

At this point we have five invariant first integrals which can be induced to the quotient. One needs one more first integral, which can be taken from Killing vector fields. 
To compute a right-invariant vector field, write the operation on the Lie group, which for coordinates $x_1,x_2,x_3,z_1,z_2,z_3$ is given by

\smallskip

$(x_1,x_2,x_3,z_1,z_2,z_3) (x_1',x_2',x_3',z_1',z_2',z_3')=(x_1+x_1', x_2+x_2', x_3+x_3',$

$ z_1+z_1'+\frac12(x_1x_2'-x_1'x_2), z_2+z_2'+\frac12(x_2x_3'-x_3x_2'),z_3+z_3'+\frac12(x_3 x_1'-x_1x_3')).$

\smallskip

 So the right-invariant vector field $V_1^*$  at $p=(x_1,x_2,x_3,z_1,z_2,z_3)$ is given by $V_1^*(p)=\partial_{x_1} + \frac12 x_2 \partial_{z_1} -\frac12 x_3 \partial_{z_3}$, which implies that $f_{V_1^*}(p,Y)=\la V_1 + x_2 Z_1 - x_3 Z_3,Y\ra$. Its gradient
$$\grad f_{V_1^*}(p,Y)=(\la Y,Z_1\ra V_2 -\la Y,Z_3\ra V_3, V_1 +  x_2 Z_1 - x_3 Z_3).$$ 
Straighforward computations show that the functions $f_{Z_i}, G, E, f_{V_1^*}$, for  i=1,2,3, are linear independent in a dense subset of $TN_{K_3}$. 

This proves that {\em  the geodesic flow on $TN_{K_3}$ is completely integrable}. But according to the non-integrable condition on $\nn_{K_3}$, the geodesic flow  cannot be Liouville integrable on $T(\Gamma\backslash N_{K_3})$ for any cocompact discrete subgroup $\Gamma < N$.

\begin{rem}  If the smoothly closed geodesics in a nilmanifold $\Gamma \backslash N$ are dense, then the nilmanifold has the density of closed geodesics property (DCG). In \cite{DDM}
the authors give conditions on the graph $G$ and on a lattice $\Gamma\subset N$ for which the quotient $\Gamma \backslash N$, a compact nilmanifold, has a dense
set of smoothly closed geodesics. 

Particularly in the situation of the graph $K_3$,  it is proved the following. Let $\Gamma$  be the lattice in $N_{K_3}$ given by $\exp(\Lambda)$ where $\Lambda$ is the vector lattice in $\nn$ given by $\Lambda = span_{2\pi \ZZ}\{ \beta \}$, for $\beta$  the orthonormal basis determined by the graph. Then  the quotient $\Gamma \backslash N_{K_3}$ has the density of closed geodesics property.
\end{rem} 

\section{Geodesic flow and  graphs} Here we study the integrability of the geodesic flow on 2-step nilpotent Lie groups  arising from  graphs. 
We consider two situations: the family of star graphs on $k+1$ vertices and  graphs in    $j$  vertices, with $j\leq 4$. We show complete integrability of geodesic flows on compact manifolds induced from star graphs so as on the corresponding simply connected Lie groups. The corresponding Lie algebras are  singular. 

For graphs with $j$ vertices, $j\leq 4$,    Liouville integrability is proved for the almost non-singular cases. 

\subsection{Star graphs}

Let $S_k$ be the star graph on $k+1$ vertices introduced in Example \ref{star}. Let $N_{S_k}$ denote the simply connected 2-step nilpotent Lie group for the Lie algebra associated to it. 
Consider its presentation  by the underlying manifold $\RR^{2k+1}$ as follows. Let $v=(x_0, x_1, \hdots,  x_k), \,  v'=(x_0', x_1', \hdots, x_k')\in \RR^{k+1}$, the group operation on $\RR^{2k+1}$ is given by
\begin{equation}\label{op-star}
\begin{array}{rcl}
(v,z_1, z_2, \hdots, z_k) (v',z_1', z_2', \hdots, z_k') & = &(x_0+x_0', x_1+x_1', \hdots, x_k+x_k', \\
& &  z_1 + z_1' +\frac{1}{2}(x_0x_1'- x_0' x_1), \\
& & z_2 + z_2' + \frac{1}{2}(x_0x_2'- x_0' x_2), \\
& & \vdots \\
& &  z_k + z_k' + \frac{1}{2}(x_0x_k'- x_0' x_k)). 
\end{array}
\end{equation}
Denote by $\partial_u$ the partial  derivation on $\RR^{2k+1}$ with respect to the variable $u$. 
A basis of left-invariant vector fields is given by 
$$
V_0(p)  =  \partial_{x_0} - \frac12 x_1 \partial_{z_1} -   \frac12 x_2 \partial_{z_2} \hdots -  \frac12 x_k \partial_{z_k} \quad 
V_i(p)  =  \partial_{x_i} + \frac12 x_0 \partial_{z_i}$$
$$
Z_i(p)  =  \partial_{z_i} \qquad \mbox{ for all }i=1, \hdots, k, 
$$  
where $p=(x_0, x_1, \hdots,  x_k, z_1, \hdots, z_k)\in N$. 
These vector fields satisfy the non-trivial Lie bracket relations
$$[V_0, V_i]= Z_i \qquad \quad \mbox{ for all }i=1, 2,\hdots, k.$$
Consider the metric on $\RR^{2k+1}$ which makes of this set an orthonormal basis. In canonical coordinates such metric  is given by$$g= (1 + \frac14 \sum_{j=1}^k x_j^2) dx_0^2 + \sum_{j=1}^k (\frac{x_0x_j}4 dx_j - \frac{x_j}2 dz_j) dx_0 + \sum_{i=1}^k(1 + \frac{x_0^2}4)dx_i^2 - \frac12 \sum_{i=1}^k x_0 dz_i dx_i + \sum_{i=1}^k dz_i^2.$$

Notice that the exponential map $\exp:\nn_{S_k} \to N_{S_k}$ is
$$\exp(\sum_{i=0}^k x_i V_i + \sum_{j=1}^k z_j Z_j) =(x_0,x_1,\hdots, x_k, z_1, z_2, \dots, z_k)$$
where $V_i, Z_j$ denote the left-invariant vector fields above.

Any right-invariant vector field on $N_{S_k}$ can be regarded as a Killing vector field.
In particular, we have the following basis of right-invariant vector fields
$$
V_0^*(p)  =  \partial_{x_0} + \frac12 x_1 \partial_{z_1} +   \frac12 x_2 \partial_{z_2} \hdots +  \frac12 x_k \partial_{z_k} \quad 
V_i^*(p)  =  \partial_{x_i} - \frac12 x_0 \partial_{z_i}$$
$$
Z_i^*(p)  =  \partial_{z_i} \qquad \mbox{ for all }i=1, \hdots, k. 
$$ Notice that $\partial_z$ is both left and right-invariant.   We induce smooth functions on $TN$ given by

\begin{equation}\label{first-int}
\begin{array}{rcl}
f_{V_0^*}(p,Y) & = &  \la V_0 + \sum_{i=1}^k x_i Z_i, Y\ra, \\
f_{V_j^*}(p,Y) & = &  \la V_j - \la W,V_0\ra Z_j, Y\ra, \\
f_{Z_j}(p,Y) & = & \la Z_j, Y\ra,
\end{array}
\end{equation}
which are first integrals of the geodesic flow, for all $j=1, \hdots, k$ and for  $exp(W)=p \in N_{S_k}$. 
It is not hard to see that for $j=1, \hdots, k$ the corresponding gradient fields are given by
$$\begin{array}{rcl}
\grad f_{V_j^*}(p,Y)   & = & (-\la Y,Z_j\ra V_0, V_j - \la W,V_0\ra Z_j)\\
\grad f_{Z_j}(p,Y)   & = & (0, Z_j)
\end{array}
$$
which are linearly independent whenever (a) $\la Y, V_0\ra \neq 0$ or (b) $\la Y, V_0\ra =0$ and $\sum_{j=1}^k \la Y, V_j\ra \la Y, Z_j\ra\neq 0$. This follows from the computations.  In fact, for the first component  we have $\sum_{i=1}^k 
a_i \la Z_i,Y\ra =0$ and on the other side, on the second component:
\begin{equation}\label{li}\sum_{i=1}^k a_i (V_i-\la W,V_0\ra Z_i) + \sum_{i=1}^k b_i Z_i + cY=0.
\end{equation}
Notice that there is only one term involving $V_0$ in Equation (\ref{li}), and the coefficient   is  $\la Y,V_0\ra$. So if $\la Y,V_0\ra\neq 0$ then $c=0$ and so $\sum_{i=1}^k a_i V_i = 0$ implies that $a_i=0$ for all $i$, and from this $b_i=0$ for all $i=1, \hdots k$.  If $\la Y,V_0\ra=0$ the  condition (b) asserts the linear independence.  This proves the first part of the following result. 

\begin{lem} Let $S_k$ denote the star graph on $k+1$ vertices. The smooth functions on $T N_{S_k}$ denoted by $f_{V_j^*}$ (as in Equation \ref{first-int}) are pairwise in involution for $j=1,\hdots k$. 

Moreover the geodesic flow on $T N_{S_k}$ is completely integrable (in the Liouville sense) since the set of first integrals $\{E, f_{Z_j}, f_{V_j^*}\}_{j=1}^k$ satisfies that any pair of first integrals is in involution and the corresponding gradients are linearly independent on an open dense set. 
\end{lem}

We only have to prove that $\{f_{V_j^*}, f_{V_i^*}\}=0$ for all $i,j=1, \hdots, k$. In fact for $i\neq j$ straighforward computations show that 
$$
\{f_{V_j^*}, f_{V_i^*}\} (p,Y) = \la Y, [V_i, V_j]\ra = 0,
$$ 
that finishes the  proof of the lemma. 

Note that all first integrals above are polynomial functions of degree one or two. In fact writing $W=\sum_{i=0}^k w_i V_i+ \sum_{i=1}^k u_i Z_i$ and $Y=\sum_{j=0}^k y_j V_j + \sum_{j=1}^k z_j Z_j$, the first integrals follow: 
$$
\begin{array}{rcl}
g(p,Y) & = & \frac12 (y_0^2+ y_1^2 + \hdots + y_k^2+z_1^2+ \hdots +z_k^2)\\
f_{V_1^*}(p,Y) & = & y_1 - w_0 z_1\\
& \vdots & \\
f_{V_k^*}(p,Y) & = & y_k - w_0 z_k\\
f_{Z_1}(p,Y) & = & z_1\\
& \vdots & \\
f_{Z_k}(p,Y) & = & z_k.
\end{array}
$$ 

\begin{rem} Take coordinates $(x_0, x_1, \hdots, x_k, z_1, \hdots, z_k)$ for $p\in N_{S_k}$ and coordinates $(y_0, y_1, \hdots, y_k, t_1, \hdots, t_k)$ on $\nn_{S_k}$ ($=\vv\oplus \zz$) relative to a basis of left-invariant vector fields. Let $F:TN_{S_k} \to \RR^{2k+1}$ be given by
$F(p,Y)=(E(p,Y), f_{Z_1}(p,Y), \hdots, f_{Z_k}(p,Y)$, $f_{V_1^*}(p,Y), \hdots, f_{V_k^*}(p,Y))$. Let $c\in \RR^{2k+1}$, namely $c=(C_0, U_1, \hdots, U_k, T_1, \hdots, T_k)$, the set $F^{-1}(c)$ gives a symplectic leaf on $TN_{S_k}$. In fact if $(p,Y)\in F^{-1}(c)$ one has
\begin{itemize}
\item $E(p,Y)=\frac12(y_0^2+ \hdots + y_k^2+ t_1^2+ \hdots t_n^2)= C_0$ and 

\item $t_i=T_i$ for all $i=1,\hdots, k$ and $y_j-x_0 T_j=Y_j$,  for $j=1,\hdots k$,
\end{itemize}
so that  coordinates $x_1, \hdots, x_k, z_1, \hdots, z_k$ have no restriction to belong to $F^{-1}(c)$. 

\end{rem}

Let $\gamma$ denote a geodesic on $N_{S_k}$. Set $\gamma(t)=exp(X(t)+Z(t))$ where $X(t)\in \vv$ and $Z(t)\in \zz$ with initial condition $X_0+Z_0$, satisfy the following system of equations (see \cite{Eb})
\begin{equation}\label{geodesic}
\begin{array}{rcl}
x_0'' & = & -a_1 x_1' - a_2 x_2' - \hdots - a_k x_k'\\
x_1'' & = & a_1 x_0' \\
 & \vdots & \\
x_k'' & = & a_k x_0' \\
z_1' & = & a_1 +\frac12(x_0 x_1'-x_0'x_1) \\
& \vdots \\
z_k' & = & a_k +\frac12(x_0 x_k'-x_0'x_k)
\end{array}
\end{equation}
where $Z_0=a_1 Z_1 + a_2 Z_2 + \hdots + a_k Z_k$ and $X(t)=\sum x_i(t) V_i$ and $Z(t)= \sum_j z_j(t) Z_j$. The map $j(Z_0)$ showed in Example \ref{star} is singular. Its kernel is 

$\{V\in \vv\,:\, V=\sum_{j=0}^k v_j V_j,\, \mbox{where } \la V, V_0\ra =0 \mbox{ and } (v_1, \hdots,v_k) \cdot (a_1, \hdots, a_k)=0 \},$ denoting  with $\cdot$  the usual inner product  for vectors in $\RR^k$. Thus if $\{g^t\}$ denotes the geodesic flow in $TN$, for every $n\in N$, and $X_0\in \vv, Z_0\in \zz$, then 
$$g^t(dL_n(X_0+Z_0)= dL_{\gamma(t)}( e^{tj(Z_0)}X_0 + Z_0),$$
where $\gamma(t)$ denotes the unique geodesic with $\gamma'(0)= dL_n(X_0+Z_0)$. 

\medskip

A Riemannian compact manifold arises as a quotient $\Lambda \backslash N_{S_k}$ where $\Lambda$ is a discrete cocompact subgroup of $N_{S_k}$. In fact $\Lambda \backslash N_{S_k}$ becomes a Riemannian manifold with the metric that makes the projection $\pi: N_{S_k} \to \Lambda \backslash N_{S_k}$ a Riemannian submersion.

Each  $2k+1$-tuple $(r,m)=(r,r_1, \ldots, r_k, m_1, \hdots, m_k) \in (\mathbb Z)^{2k+1}$ defines a lattice in $N_{S_k}$ by
\begin{equation}\label{eq:lattice}
  \Lambda_{(r,m)} = rm_0 \ZZ \times 2r_1 \ZZ \times  \hdots \times 2r_k \ZZ \times m_1 \ZZ \times m_2 \ZZ \times \hdots \times m_k\ZZ, 
\end{equation}
for $m_0= m_1 m_2 \dots m_k$.

 Note that there are non-isomorphic lattices in this family so that we get many non-diffeomorphic compact manifolds.

Since the quotient projection $\pi: N_{S_k} \to \Lambda_{(r,m)} \backslash N_{S_k}$ is a Riemannian submersion and furthermore a local isometry, we can identify the tangent bundle  $T (\Lambda_{(r,m)} \backslash N_{S_k})$ with $(\Lambda_{(r,m)} \backslash N_{S_k}) \times \mathfrak \nn_{S_k}$. The projection $\pi$ maps geodesics into geodesics and  the energy function $\tilde E$ on $T(\Lambda_{(r,m)} \backslash N_{S_k})$ is related to the energy function $E$ on $T N_{S_k}$ by
$$\tilde E(\Lambda_{(r,m)} p, Y) = E(p, Y) = \frac{1}{2}\langle Y, Y\rangle$$
and clearly it is well defined. 

All invariant first integrals on $T N_{S_k}$ descend to the quotients, since they do not depend on  the  coordinates of $p\in N_{S_k}$. One defines
$$\tilde f_{Z_j}( \Lambda_{(r,m)}p, Y) = f_{Z_j}(p, Y), \qquad \mbox{ for all }j=1, \dots, k,$$
which are first integrals of the geodesic flow of $T(\Lambda_{(r,m)} \backslash N_{S_k})$. Moreover, such first integrals are in involution, since for all $f, g \in C^\infty(T(\Lambda_{(r,m)} \backslash N_{S_k}))$ we have
$$\{f \circ \pi, g \circ \pi\} = \{f, g\} \circ \pi.$$

Note that the integrals $f_{V_j^*}$, $j = 1, \ldots, k$  do not descend directly to the quotient. However one can construct first integrals on the quotient with the following argument. Let $(p, Y) \in T N_{S_k}$ and $q \in \Lambda_{(r,m)}$. 
Take $W, W' \in \mathfrak n_{S_k}$ such that $\exp W = p$, $\exp W' = q$. Observe that $(W+W')_{\vv}= W_{\vv} + W'_{\vv}$, where $U_{\vv}$ denotes the orthogonal projection of $U\in \nn_{S_k}$ over $\vv=\zz^{\perp}$.  So we get
\begin{align*}
  f_{V_j^*}(qp, Y) & = \langle Y, V_j\rangle - \langle Z_j, Y\rangle\langle (W  +  W'), V_0\rangle \\
  & = f_{V_j^*}(p, Y) - f_{Z_j}(p, Y)\langle W', V_0\rangle.
\end{align*}
Since $\langle W', V_0\rangle \in \mathbb Z$ we have that 
$$f_{V_j^*}(qp, Y) = f_{V_j^*}(p, Y) \mod f_{Z_j}(p, Y)\mathbb Z$$ for every $j=1, \dots, k$
and since $f_{Z_j}$ is a first integral of the geodesic flow, we have that the function
$$\hat f_{V_j^*} (p, Y) = \sin\left(2\pi\frac{f_{V_j^*}(p, Y)}{f_{Z_j}(p, Y)}\right)$$
descends to $\Lambda_{(r,m)} \backslash N_{S_k}$ and is constant along the integral curves of the geodesic vector field in $T(\Lambda_{(r,m)} \backslash N_{S_k})$. In order to get a smooth first integral let 
$$\bar F_j(p, Y) = e^{-1/f_{Z_j}(p, Y)^2}\hat f_j(p, Y)$$
and let us define
\begin{equation}\label{fcompact}
\tilde{F}_j(\Lambda_rp, Y) = \bar F_j(p, Y).
\end{equation}

So the functions $\tilde{F}_k$ are smooth (non-analytic) first integrals for the geodesic flow on $T(\Lambda_{(r,m)} \backslash N_{S_k})$. It follows from a direct calculation making use of properties of the Poisson bracket  that the families $f_{Z_i},  \tilde{F}_j$, $i= 1, \ldots, k$ are in involution. In fact for a pair of differentiable functions on $M$ and for $h:\RR \to \RR$, it holds $\{f, h\circ g\}=h'\{f,g\}$. To prove the linear independence, notice that the new gradients can be written in terms of the gradients on $N_{S_k}$, which are multiplied by differentiable real functions. The independence follows by asking the corresponding determinant is not trivial on the right open set. So the geodesic flow in $T(\Lambda_{(r,m)} \backslash N_{S_k})$ is completely integrable in the sense of Liouville. 

\begin{thm}
Let $N_{S_k}$ be the 2-step nilpotent Lie group attached to the star graph in $k+1$ vertices $S_k$,  endowed with the standard metric and let $\Lambda_{(r,m)}$ denote the lattice  in (\ref{eq:lattice}). If $\Lambda_{(r,m)} \backslash N_{S_k}$ is the corresponding compact manifold with the induced metric, then the geodesic flow in $T(\Lambda_{(r,m)} \backslash N_{S_k})$ is completely integrable with smooth first integrals $\{E, f_{Z_i}, \tilde{F}_i\}$, for $i=1, \hdots, k$.
\end{thm}

See  explanations on the proof and the  topology of these compact manifolds in \cite{Bu2}, where the author worked at the cotangent. Another difference is the presentation of the Lie group. See the isomorphism in Equation \ref{psi} below. 

\smallskip

Let  $\tilde{H}\subset N_{S_k}$ be the normal subgroup of dimension $2k$ defined as $\tilde{H}=\{g\in N_{S_k} \, : \, g=(v,z) \quad \mbox{and}\quad v=(0,x_1, \hdots, x_k)\}$. Note that $\tilde{H}$ is  abelian  and $\tilde{\Lambda}_{(r,m)} = \tilde{H}\cap \Lambda_{(r,m)}$ is a lattice in $\tilde{H}$. So $\tilde{\Lambda}_{(r,m)} \backslash \tilde{H} \simeq T^{2k}$.

\medskip

Note that $N_{S_1}$ is (isomorphic to) the Heisenberg Lie group $H_3$. A known  presentation of $H_3$ is given in terms of matrices as the set:
 $$\left\{ \left( 
\begin{matrix}
1 & x & z \\
0 & 1 & y \\
0 & 0 & 1 
\end{matrix} \right)
 \mbox{ for }x,y, z\in \RR
\right\}$$
together with the usual product of matrices. 
The subgroup  $\Gamma_r$ consisting of matrices of the form
$$\Gamma_r = \left\{ \left( 
\begin{matrix}
1 & rn & q \\
0 & 1 & m \\
0 & 0 & 1 
\end{matrix} \right)
 \mbox{ for } m, n, q\in \ZZ
\right\}$$
for a fixed $r\in \NN$ gives rice to a cocompact lattice in $H_3$. The lattice $\Gamma_r$ induces an action on $H_3$ so that the class of $(x,y,z)\in H_3$ is
$\overline{(x,y,z)}=\{(x+rn, y+m, z+rn + s)\, :\, n\in \ZZ, m\in \ZZ, s\in \ZZ\}$. 

Denoting by $\bar{\Gamma}_r$ the subgroup isomorphic to $\ZZ$  and  by $\tilde{\Gamma}$ the subgroup isomorphic to $\ZZ^2$ given respectively by
$$\bar{\Gamma}_r=  \{(rn, 0, 0):n\in \ZZ\} \simeq r \ZZ, \qquad \tilde{\Gamma}=\{(0,m,q): m, q\in \ZZ\}\simeq \ZZ \times \ZZ$$
we get the semidirect product group $\bar{\Gamma}_r \ltimes \tilde{\Gamma}$, where the action is given by
$ rn \cdot (m,q)=(m, q + rnm)$, via identifications. 
So the map $\Psi : \bar{\Gamma}_r \ltimes \tilde{\Gamma} \to \Gamma_r$ given as $(rn, (m,q)) \to (rn, m,q)$ is a group isomorphism. And the action of $\Gamma_r$ on $H_3$ by translations on the left is equivalent to an action of $\bar{\Gamma}_r \ltimes \tilde{\Gamma}$ on $H_3$. 

Now the action of $\tilde{\Gamma}$ on $H_3$ gives
$$(0,m, q) \cdot (x,y,z)=(x, y+m, z+q) \quad \mbox{ so that } \tilde{\Gamma} \backslash H_3 \simeq \RR \times T^2.$$ 
One also has
$$(rn, 0,0) \cdot (0,m, q) \cdot (x,y,z) = (0,m, q) \cdot (rn,0,0) \cdot (x,y,z).$$ 
Finally the action of  $\bar{\Gamma}_r$ on $\RR \times T^2$ induces
the action of $S^1$ on $\RR\times T^2$ and  $\Gamma_r\backslash H_3$  is a $S^1$-fiber bundle over $T^2$:
$$S^1 \quad \to \quad \Gamma_r\backslash H_3 \quad \to \quad T^2.$$

A similar procedure generalizes to $H_{2n+1}$ showing  that Heisenberg nilmanifolds $\Gamma\backslash H_{2n+1}$, as topological spaces  are $T^{n}$-fiber bundles over $T^{n+1}$ see \cite{BJ}. 

Analogously one shows  the fibration we get from the nilmanifolds arising from every star graph. In fact, one defines similar subgroups $\Gamma_r$, $\bar{\Gamma}_r$ and $\tilde{\Gamma}$ in $N_{S_k}$. This is  explained above. 

A presentation of $N_{S_k}$ is given by the $(k+2)\times(k+2)$-matrices of the form
$$\left\{ 
\left( \begin{matrix}
1 & x_0 & z_1 & z_2 & \hdots & z_k \\
0 & 1 & x_1 & x_2 & \hdots & x_k \\
0 & 0 & 1  & 0 & \hdots & 0 \\
0 & 0 & 0 & 1 & \hdots   & 0 \\
0 &  &  & 0 & \ddots   & \vdots \\
0 &  & & &  & 1
\end{matrix}
\right),
\right\}
$$
where we think in the usual matrix multiplication. 
The map 
\begin{equation}\label{psi}
\Psi: \left( \begin{matrix}
1 & x_0 & z_1 & z_2 & \hdots & z_k \\
0 & 1 & x_1 & x_2 & \hdots & x_k \\
0 & 0 & 1  & 0 & \hdots & 0 \\
0 & 0 & 0 & 1 & \hdots   & 0 \\
0 &  &  & 0 & \ddots   & \vdots \\
0 &  & & &  & 1
\end{matrix}
\right)
 \quad \mapsto \quad (x_0,x_1, x_2, \hdots, x_k, z_1-\frac12 x_0x_1,  \hdots, z_k-\frac12 x_0x_k)
\end{equation}
gives an isomorphism between the above group and that one defined with the multiplication operation of Equation \ref{op-star}.

For a fixed $r$, take  the lattice $\Gamma_r\subset N_{S_k}$ given by matrices of the form
$$\Gamma_r= \left\{ 
\left( \begin{matrix}
1 & rn & q_1 & q_2 & \hdots & q_k \\
0 & 1 & m_1 & m_2 & \hdots & m_k \\
0 & 0 & 1  & 0 & \hdots & 0 \\
0 & 0 & 0 & 1 & \hdots   & 0 \\
0 &  &  & 0 & \ddots   & \vdots \\
0 &  & & &  & 1
\end{matrix}
\right) \qquad \mbox{ for } \quad n, m_i, q_j\in \ZZ, \forall i,j\right\}.
$$
One starts with  the abelian subgroup $\tilde{\Gamma}=\{ (0,m_1, m_2, \hdots, m_k, q_1, q_2, \hdots, q_k) : m_i, q_j\in \ZZ\}$, which acts on the left on  $N_{S_k}$ by
$$\begin{array}{rcl}
(0,m_1, m_2, \hdots, m_k, q_1, q_2, \hdots, q_k)\cdot (x_0, x_1, \hdots, x_k, z_1, \hdots, z_k) & = & (x_0, x_1+m_1, \hdots,\\
& &  x_k+m_k, z_1+q_1, \\
& & \hdots, z_k+q_k).
\end{array}
$$
This shows that $\tilde{\Gamma} \backslash N_{S_k}\simeq \RR \times T^{2n}$. Let $X_0\in N_{S_k}$ denote the element $X_0=(x_0, x_1+m_1, \hdots, x_k+m_k,z_1+q_1, \hdots, z_k+q_k)$ and  take the action of $\bar{\Gamma}_r=\{(rn,0, \hdots, 0):n\in \ZZ\}< N_{S_k}$: 
\begin{equation}\label{result}
\begin{array}{rcl}
(rn, 0, \hdots, 0)\cdot X_0 & = & (x_0+rn, x_1+m_1, \hdots, x_k+m_k, \\
& &  z_1+q_1+rnx_1+rnm_1, \hdots,  \hdots,  z_k+q_k+rnx_k+rnm_k).
\end {array}
\end{equation}
 
As above for $H_3$ one has ${\Gamma}_r \simeq \bar{\Gamma}_r \times \tilde{\Gamma}$ and the action of $\Gamma_r$  on $N_{S_k}$ translates into an action of $\bar{\Gamma}_r \ltimes \tilde{\Gamma}$ on $N_{S-k}$ showing that 
 the compact manifold $\Gamma_r \backslash N_{S_k}$ is a $S^1$-fiber bundle over $T^{2k}$:
$$S^1 \quad \to \quad  \Gamma_r \backslash N_{S_k} \quad \to \quad T^{2n}.$$
The action of $S^1$ on each torus $T^2$ is  induced by the matrix action on each subspace $(z_i,x_i)$ for $i=1, \hdots, k$:
$$\left(
\begin{matrix}
1 & r \\
0 & 1 
\end{matrix}
\right), 
$$
as already explained on the Heisenberg Lie group $H_3$. 

\begin{rem} Let $\varphi_t$ denote the geodesic flow. It was proved in \cite{Bu2} that $\varphi_t$ is non-degenerate in the sense of KAM theory. Moreover the topological entropy of $\varphi_t$ vanishes. Indeed $\pi_1( \Gamma_1 \backslash N_{S_k})$ has no abelian subgroup of finite index. For more information on the topology see \cite{Bu2}. 
\end{rem}

\subsection{Non-commutative integrability and graphs} Here we prove the integrability of geodesic flows on manifolds associated to graphs in a low number of vertices, $k$ where $k\leq 4$. 

\smallskip

Let $H: TN \to \RR$ denote a smooth function. One says that it is {\em integrable in the non-commutative sense of Nekhorosev} or simply {\em integrable} if one has the following conditions. Assume $F = (H = f_1, \hdots ,f_{n-k}, g_1, \hdots ,g_{2k})$ is a smooth
map on $TN$ where $\dim N=n$ and $k\geq  0$, that satisfies the three conditions:
\begin{enumerate}[(i)]
\item $ rank\, dF = n + k$ on an open, dense subset of $TN$;
\item for all $a, b = 1, \hdots ,n - k$ and all $c = 1, \hdots, 2k$ : $\{f_a, f_b\} = \{f_a, g_c\} = 0$;
\item for each regular value $c \in \RR^{n+k}$, each connected component of $F^{-1}(c)$ is compact. 
\end{enumerate}

In this situation the Theorem of Nekhorosev describes the level sets of $F$, and the flow for the Hamiltonian $X_H$ in terms of a flow on the torus. In \cite{Bu1} Butler proves integrability in the non-commutative sense of the geodesic flow for $D \backslash N$, where $D$ is a cocompact lattice on $N$ and $N$ is 2-step nilpotent Lie group whose Lie algebra is almost non-singular. In fact one can find $3s+t$ first integrals, where $\dim \vv=2s$ and $\dim \zz=t$. These functions can be obtained by the two ways studied in this section, those invariant functions from  Proposition \ref{fin} and the $n$ functions arising from a basis a right-invariant vector fields. 

Since Bolsinov and Jovanovic \cite{BJ} proved that integrability in the non-commutative sense implies Liouville integrability, the previous results of Butler give the Liouville integrability for an important family of 2-step nilpotent Lie groups and their compact quotients.

Let us explain the construction. Let $N$ denote a Lie group equipped with a left-invariant metric and with Lie algebra $\nn=\vv\oplus \zz$ with $\dim \vv=2n$, $\dim \zz=m$. Making use of Killing vectors we construct $2n+m$ first-integrals of the geodesic low. With Proposition \ref{fin} we construct $n$  invariant first integrals. Assume this is a Lie algebra. In this Lie algebra we have $n+m$ which are in involution. So that this gives:
$2n+m + n + n + m= 2m + 4n=\dim TN$. This gives the complete integrability for almost non-singular Lie algebras.  As corollary one gets the  next result.

\begin{cor} Let $G$ denote a connected graph on $k$ vertices with $k\leq 4$. Then except for the complete graph $K_3$, any 2-step nilpotent Lie group $N_G$ so as the corresponding compact quotient admits  a completely integrable geodesic flow. 
\end{cor}

 Take the graphs in Example \ref{graphsandvertex}: 
For two vertices, we have the Heisenberg Lie algebra of dimension three whose geodesic flow is completey integrable with any left-invariant metric. 
For three vertices, the connected graphs are $S_3$ and $K_3$ which were explained above. We need to concentrate in algebras coming from graphs with four vertices. 

\begin{example} Let $P$ the path in four vertices. Let $N_{P}$ denote the 2-step nilpotent Lie group associated to $P$. Then $N_P$ has dimension seven. Take $p=exp(W)$ and the functions on $TN_P$ given by
$$\begin{array}{rcl}
f_{V_1^*}(p,Y)  & = &    \la V_1 + \la W,V_2\ra Z_1, Y\ra, \\
f_{V_2^*}(p,Y)  & = &    \la V_2 - \la W,V_1\ra Z_1 + \la W, V_3 \ra Z_3, Y\ra, \\
f_{V_3^*}(p,Y)  & = &    \la V_3 - \la W,V_2\ra Z_2 + \la W, V_4\ra Z_3, Y\ra, \\
f_{V_4^*}(p,Y)  & = &  \la V_4 - \la W,V_3\ra Z_3, Y\ra.
\end{array}
$$
In this situation we do not need the functions of Proposition \ref{fin}. In fact Killing vectors obtained as right-invariant vector fields give rise to a subalgebra of first integrals of dimension 7. Among these first integrals, we have 3 from the center which are in involution and also two more functions, $f_{V_1^*}$ and $f_{V_4^*}$ which are in involution. In this situation since $\nn_P$ is almost non-singular, we could apply the result in \cite{BJ} to get the complete integrability. But we are also able to give explicitly a family of first integrals in involution. In fact, we already have five first integrals in involution as in Example \ref{path}. Choose the first integrals above $f_{V_1^*}, f_{V_4^*}$. It is not hard to prove that the corresponding gradient fields are given by:

\smallskip

$\grad f_{V_1^*}(p,Y) = (\la Y, Z_1\ra V_2, V_1+\la W,V_2,V_2\ra Z_1)$

\smallskip

$\grad f_{V_4^*}(p,Y) = (-\la Y, Z_3\ra V_3, V_4 - \la W,V_3\ra Z_3).$

\smallskip

Making use of this information and by doing similar computations as already shown for Lie groups associated to star graphs, one can prove the following:

\smallskip

{\em The set $E,h,f_{Z_1}, f_{Z_2}, f_{Z_3}, f_{V_1^*},f_{V_4^*}$ is a set of first integrals in involution. } 

\smallskip The map $h$ above were defined in Example \ref{path}.
\end{example}

Acknowledgement: The author thanks the generous comments and suggestions of the referee, to improve  the first version of the work.


\begin{thebibliography}{GGGG}




\bibitem{BJ}{\sc A. Bolsinov,  B. Jovanovi\'c}, { Noncommutative Integrability, Moment Map
and Geodesic Flows}, Annals of Global Analysis and Geometry {\bf 23}, (2003) 305--322.


\bibitem{BJ2} {\sc A.  Bolsinov, B. Jovanovi\'c}, Complete involutive algebras of functions on cotangent bundles of homogeneous spaces, Math. Z. {\bf 246}, (2004) 213--236.

\bibitem{BT1} {\sc A. Bolsinov, I.  Taimanov}, Integrable geodesic flows with positive topological entropy, Invent. Math. {\bf 140}, (2000) 639--650.

\bibitem{BT2} {\sc A. Bolsinov,I. Taimanov}, On an example of an integrable geodesic flow with positive topological entropy, Russ. Math. Surv. {\bf 54}(4), (1999) 833--835.

\bibitem{Bu1} {\sc L. Butler}, { Integrable geodesic flows with wild first integrals: the case of two-step 
nilmanifolds}. Ergodic Theory Dynam. Systems {\bf 23} (3), (2003) 771--797.

\bibitem{Bu2} {\sc L. Butler}, New examples of integrable geodesic flows, Asian J. Math. 
{\bf 4} (3), (2000) 515--526. 

\bibitem{Bu3} {\sc L. Butler}, Zero entropy, non-integrable geodesic flows and a non-commutative rotation vector, Transactions of the AMS {\bf 355} (9), (2003) 3641--3650.

\bibitem{DM} {\sc S. G. Dani, M. G. Mainkar}, Anosov automorphisms on compact nilmanifolds associated with graphs, Trans. Amer. Math. Soc. {\bf 357},  (2005) 2235--2251.

\bibitem{DDM} {\sc R.  DeCoste, L. Demeyer, M. Mainkar}, Graphs and metric 2-step nilpotent Lie algebras, Adv. Geom. {\bf 18} (3), (2018) 265--284 .

\bibitem{Eb} {\sc P. Eberlein}, {Geometry of 2-step nilpotent Lie groups with a left-invariant metric}, Ann. Sci. E. N. S., 4 serie, {\bf 27} (5), (1994) 611--660. 
\bibitem{Eb2} {\sc P. Eberlein}, {Left invariant geometry of Lie groups}, Cubo 6 (1), (2004) 427--510

\bibitem{GMS} {\sc  C. Gordon, Y. Mao, D. Schueth}, { Symplectic rigidity of geodesic flows on two-step nilmanifolds }, Ann. Scient. \'Ecole Norm. Sup. (4) {\bf 30} (4),  (1997) 417--427.

\bibitem{GM} {\sc R. Gornet and M. Mast}, The length spectrum of riemannian two-step nilmanifolds, Ann. Scient. \'Ecole Norm. Sup. (4) {\bf 33} (2), (2000) 181--209.

\bibitem{KOR} {\sc A. Kocsard, G. P. Ovando, S. Reggiani}, On first integrals of the geodesic flow on Heisenberg nilmanifolds, Diff. Geom. Appl. {\bf 49}, (2016) 496--509. 


\bibitem{Ko2}{\sc B. Kostant}, { The solution to a generalized  Toda lattice
and representation theory}, Advances in Math. {\bf 39},  (1979)195--338.

\bibitem{Ma} {\sc M. Mainkar}, Graphs and two-step nilpotent Lie algebras, Groups Geom. Dyn. {\bf 9} (1), (2015) 55--65.

\bibitem{Mal} {\sc A. I. Mal'cev}, {\it On a class of homogeneous spaces}, Amer. Math. Soc.Translations {\bf  39}, 1951.

\bibitem{Pa} {\sc G. Paternain, R.J.  Spatzier},  New examples of manifolds with completely integrable flows, Adv. Math. {\bf 108}, (1994) 346--366.


\bibitem{Sc} {\sc D. Schueth}, Integrability of geodesic flows and isospectrality of Riemannian manifolds, Math. Z. {\bf 260} (3), (2008) 595--613 . 

\bibitem{Sy}{\sc W. Symes}, {Systems of Toda type, inverse
spectral problems and representation theory}, Invent. Math. {\bf 59} (1978), 13--53.

\bibitem{Ta1} {\sc I. A. Taimanov}, Topological obstructions to integrability of geodesic flows on non-simply-connected manifolds, Izv. Akad. Nauk SSSR, Ser. Mat. {\bf 51}(2) (1987,) 429--435.

\bibitem{Ta2} {\sc I. A. Taimanov}, Topology of Riemannian manifolds with integrable geodesic flows, Tr. Mat. Inst. Steklova {\bf 205} (1994), 150--163.

\bibitem{Th} {\sc A. Thimm}, { Integrable geodesic flows on homogeneous spaces}, Ergod. Th.  \& Dynam. Sys. {\bf 1} (1981), 495--517.

\end{thebibliography}
\end{document}